\documentclass{amsart}
\usepackage{pdfpages}
\newtheorem{theorem}{Theorem}[section]

\theoremstyle{definition}
\newtheorem{definition}[theorem]{Definition}

\theoremstyle{remark}

\numberwithin{equation}{section}

\begin{document}

\title
 {Homothetical surfaces in three 
 dimensional 
 pseudo-Galilean spaces 
 satisfying $\Delta^{II}\textbf{x}_i=\lambda_i\textbf{x}_i$}

\author[M.S. Lone]{Mohamd Saleem Lone}

\address{%
International Centre for Theoretical Sciences, \\
Tata Institute of Fundamental Research,\\
 560089, Bengaluru, India.}

\email{saleemraja2008@gmail.com, mohamdsaleem.lone@icts.res.in}





\begin{abstract}
A homothetical surface arises as a
graph of a function $z =  \varphi_1(v_1) \varphi_2(v_2)$. In this paper, we study the homothetical surfaces in three dimensional pseudo-Galilean space$\left(\mathbb{G}_3^1\right)$ satisfying the conditions $\Delta^{II}\textbf{x}_i=\lambda_i\textbf{x}_i,$ where $\Delta^{II}$ is the Laplacian  with respect to second fundamental form. In particular, we show the non-existence of any such type of surface in $\mathbb{G}_3^1.$
\end{abstract}

\maketitle
\section{Introduction}
An Euclidean submanifold is said to be of finite type (or finite Chen-type) if its coordinate function is the finite sum of eigenfunctions of its Laplacian. B.-Y. Chen posed the problem of classifying the finite type of surfaces in 3-dimensional Euclidean space $\mathbb{E}^3.$ The notion of finite type can be extended to any smooth function on a submanifold of a Euclidean space or any ambient space.\newline
  Let $\textbf{x}:\textbf{M}\rightarrow\mathbb{E}^m$ be an isometric immersion of a connected $\textit{n-}$dimensional manifold in the $\textit{m-}$dimensional Euclidean space $\mathbb{E}^m.$ Denote by $\textbf{H}$ and $\Delta$ the mean curvature and the Laplacian of $\textbf{M}$ with respect to the Riemannian metric on $\textbf{M}$ induced from that of $\mathbb{E}^m,$ respectively. Takahashi \cite{Gr15} proved that the submanifold in $\mathbb{E}^m$ satisfying $\Delta \textbf{x}=\lambda \textbf{x},$ i.e., all the coordinate functions are eigenfunctions of the Laplacian with the same eigenvalue $\lambda \in \mathbb{R}$, are either the minimal submanifolds of $\mathbb{E}^m$ or the minimal submanifolds of hypersphere $\mathbb{S}^{m-1}$ in $\mathbb{E}^m$. \newline
  As an extension of Takahashi theorem, in \cite{Gr11} Garay studied hypersurfaces in $\mathbb{E}^m$ whose coordinate functions are eigenfunctions of the Laplacian, but not necessary according to the same eigenvalue. He considered hypersurfaces in $\mathbb{E}^m$ satisfying the condition
$$\Delta \textbf{x}=\textbf{Ax},$$
where $\textit{\bf A}\in  Mat(m, \mathbb{R})$ is an $m\times m-$diagonal matrix and proved that such hypersurfaces are minimal in $\mathbb{E}^m$ and open pieces of either round hypersurfaces or generalized right spherical cylinders. \newline
  Related to this, Dillen, Pas and Verstraelen \cite{Gr10} investigated surfaces in $\mathbb{E}^3$ whose immersions satisfy the condition
$$\Delta \textbf{x}=\textbf{Ax}+\textbf{B},$$
where $\textbf{B}\in\mathbb{R}^3$. For the Lorentzian version of surfaces satisfying $\Delta \textbf{x}=\textbf{Ax}+\textbf{B}$, Alias, Ferrrandez and Lucas \cite{Gr1} proved that the only such surfaces are minimal surfaces and open pieces of Lorentz circular cylinders, hyperbolic cylinders, Lorentz hyperbolic cylinders, hyperbolic spaces or pseudo-spheres.\newline
  The notion of an isometric immersion $\textbf{x}$ is naturally extendable to smooth functions on submanifolds of Euclidean space or pseudo-Euclidean space. The most natural one of them is the Gauss map of the manifold. In particular, if the submanifold is a hypersurface, the Gauss map can be identified with the unit normal vector field to it. Baikoussis and Verstraelen \cite{Gr3} studied the helicoidal surfaces in $\mathbb{E}^3$. Choi \cite{Gr8} completely classified the surface of revolution in the three-dimensional Minkowski space $\mathbb{E}_1^3$ satisfying the condition
$$\Delta\textbf{G}=\textbf{AG}.$$
    Yoon \cite{Gr16} classified the translation surfaces in the three-dimensional Galilean space under the condition
$$\Delta\textbf{x}^i=\lambda^i{\bf x}^i,$$
\noindent where $\lambda^i\in \mathbb{R}.$ The authors in \cite{Gr4,Gr12} classified translation surface and surface of revolution, respectively in three-dimensional spaces satisfying
$$\Delta^{III}\textbf{r}_i=\mu_i{\bf r}_i.$$
 Yu and Liu \cite{Gr17} and the authors in \cite{Gr14} studied the homothetical minimal surfaces in 3-dimensional Euclidean and Minkowski spaces. Bekkar and Senoussi \cite{Gr5} classified the homothetical surfaces in 3-dimensional Euclidean and Lorentzian spaces satisfying
$$\Delta \textbf{r}_i=\lambda_i {\bf r}_i.$$
  Aydin, \"{O}\v{g}renmi\c{s} and Erg\"{u}t \cite{Gr2} investigated the homothetical surfaces in pseudo-Galilean space with null Gaussian and mean curvature. Karacan, Yoon and Bukcu \cite{Gr13} classified translation surfaces of first type satisfying $\Delta^J\textbf{x}_i=\lambda_i\textbf{x}_i$, $J=1,2$ and $\Delta^{III}\textbf{x}_i=\lambda_i\textbf{x}_i.$ Recently, Cakmak et al. \cite{Gr6} studied the translation surfaces in the three-dimensional Galilean space satisfying
$$\Delta^{II}\textbf{x}_i=\lambda_i{\bf x}_i.$$
  Motivated by all of the above research and in particular by the discussion in \cite{Gr6}, the focus of this paper is to investigate a homothetical surface in pseudo-Galilean space satisfying $\Delta^{II}\textbf{x}_i=\lambda_i{\bf x}_i$. The importance of the paper lies in the fact: there are no such surfaces with non-trivial Gaussian curvature.
\section{Preliminaries}
The pseudo-Galilean space $\mathbb{G}_3^1$ is a Cayley-Klein space defined from a three-dimensional projective space $P\mathbb{R}^3$ with the absolute figure that consists of an ordered triplet $\{\omega,f,I\}$, where $\omega$ is the ideal(absolute) plane, $f$ the line (absolute line) in $\omega$ and $I$ the fixed hyperbolic involution of the points of $f$. We introduce homogeneous coordinates in $\mathbb{G}_3^1$ in such a way that the absolute plane $\omega$ is given by $x_0=0,$ the absolute line $f$ by $x_0=x_1=0$ and the hyperbolic involution by $(0:0:x_2:x_3)\rightarrow(0:0:x_3:x_2)$. In affine coordinates defined by $(x_0:x_1:x_2:x_3)\rightarrow(1:x:y:z),$ distance between points $Q_i=(x_i,y_i,z_i),i=1,2$ is defined by: (\cite{Divjak,Sipus})
\begin{eqnarray}\label{G1}
d(Q_1,Q_2)=\left\{
             \begin{array}{ll}
              \quad\quad \quad |x_2-x_1|, \quad \quad\quad x_1\neq x_2, \\
               \sqrt{(y_2-y_1)^2+(z_2-z_1)^2}, \quad x_1=x_2.
             \end{array}
           \right.
\end{eqnarray}
The group of motions of $\mathbb{G}_3^1$ is a six parameter group given (in affine coordinates) by
\begin{eqnarray*}
&&\overline{x}=a+x,\\
&&\overline{y}=b+cx+y\cosh\theta+z\sinh \theta,\\
&& \overline{z}=d+ex+y\sinh \theta +z \cosh \theta.
\end{eqnarray*}
The pseudo-Galilean scalar product of two vectors $Q_1=(x_1,x_2,x_3)$ and $Q_2=(y_1,y_2,y_3)$ is defined as
\begin{eqnarray*}
  Q_1\cdot Q_2=\left\{
                 \begin{array}{ll}
                   x_1y_1\quad \quad \quad \quad  \textrm{ if } x_1\neq0 \textrm{ or } y_1\neq 0, \\
                   x_2y_2-x_3y_3 \quad \quad \textrm{ if } x_1=0 \textrm{ and } y_1= 0.
                 \end{array}
               \right.
\end{eqnarray*}
In pseudo-Galilean space a vector $Q=(x_1,x_2,x_3)$ is called isotropic (non-isotropic) if $x_1=0(x_1 \neq 0).$ All unit non-isotropic vectors are of the form $(1,x_2,x_3).$ The isotropic vector $Q=(0,x_2,x_3)$ is called \textit{spacelike, timelike} and \textit{lightlike} if $x_2^2-x_3^2>0$, $x_2^2-x_3^2<0$ and $x_2=\pm x_3$, respectively. The pseudo-Galilean cross product of $Q_1$ and $Q_2$ in $\mathbb{G}_3^1$ is given by
\begin{eqnarray*}
Q_1 \times Q_2=\left|
  \begin{array}{ccc}
    0 & -e_2 & e_3 \\
    x_1 & x_2 & x_3 \\
    y_1 & y_2 & y_3 \\
  \end{array}
\right|,
\end{eqnarray*}where $e_2$ and $e_3$ are the standard basis.\newline
  Let $\textbf{M},$ be a $C^r,$ $r\geq 1$ surface in pseudo-Galilean space $\mathbb{G}_3^1$ parameterized by
 $$\textbf{x}(v_1,v_2)=(x(v_1,v_2),y(v_1,v_2),z(v_1,v_2)).$$ From now onward set $x_i=\frac{\partial x}{\partial v_i},$ $\{i=1,2\}$, similarly for $y(v_1,v_2)$ and $z(v_1,v_2)$. The surface $\textbf{M}$ has the following first fundamental form
$$\textbf{I}=\left(
               \begin{array}{cc}
                 ds_1^2 & 0 \\
                 0 & ds_2^2 \\
               \end{array}
             \right),
   $$
with $ds^2=(g_1 dv_1^2+g_2dv_2^2)^2+(h_{11}dv_1^2+2h_{12}dv_1dv_2+h_{22}dv_2^2), $ $g_i=x_i$ and $h_{ij}=\widetilde{\textbf{x}_i}\cdot \widetilde{\textbf{x}_j}$ stands for derivatives of the first coordinate function $x(v_1,v_2)$ with respect to $v_1,v_2$ and for the Euclidean scalar product of the projections $\widetilde{\bf x}_k$ of the vectors $\textbf{x}_k$ onto the $yz-$plane, respectively.
A surface is called admissible if it has no Euclidean tangent planes. Therefore, for an admissible surface either $g_1\neq 0,$ or $g_2\neq 0,$ holds. An admissible surface can always be expressed as $$z= \varphi(v_1,v_2).$$
The vector $\textbf{N}$ defines a normal vector to the surface and is given by
$$\textbf{N}=\frac{1}{W}(0,-x_2z_1+x_1z_2,x_1y_2-x_2y_1),$$
where $W=\sqrt{|(x_1y_2-x_2y_1)^2-(x_1z_2-x_2z_1)^2|}$ and $\textbf{N}\cdot \textbf{N}=\epsilon=\pm1. $\newline
Hence two types of admissible surfaces can be distinguished: spacelike having timelike unit normal $(\epsilon =-1)$ and timelike having spacelike unit normal $(\epsilon =1).$
The Gaussian $\textbf{K}$ and the mean curvature $\textbf{H}$ are $C^{r-2}(r\geq 2)$ functions, defined by
$$\textbf{K}=-\epsilon\frac{LN-M^2}{W^2}, \quad \textbf{H}=-\epsilon\frac{g_2^2L-2g_1g_2M+g_1^2N}{2W^2},$$
where
$$L_{ij}=\epsilon\frac{x_1\tilde{\textbf{x}}_{ij}-x_{ij}\tilde{\textbf{x}}_1}{x_1}\cdot\textbf{N},\quad x_1=g_1\neq0.$$
For convenience, we will use $L,$ $M,$ $N$ instead of $L_{ij},$ $i,j=1,2$.\newline
  It is well known in terms of local coordinates $\{v_1,v_2\}$ of $\textbf{M}$, the Laplacian operator $\Delta^{II}$ with respect to the second fundamental form on $\textbf{M}$ is defined by \cite{Neill}
\begin{eqnarray}\label{G2}
\Delta^{II}{\bf x}=\frac{-1}{\sqrt{LN-M^2}}\left[\frac{\partial}{\partial v_1}\left(\frac{N\textbf{x}_1-M\textbf{x}_2}{\sqrt{LN-M^2}}\right)
-\frac{\partial}{\partial v_2}\left(\frac{M\textbf{x}_1-L\textbf{x}_2}{\sqrt{LN-M^2}}\right)\right],
\end{eqnarray}
where the second fundamental form is non-degenerate, or $LN-M^2\neq 0.$ 
\section{Homothetical surfaces in $\mathbb{G}_3^1$}
A surface $\textbf{M}$ in the pseudo-Galilean space $\mathbb{G}_3^1$ is called a homothetical (or factorable) surface if it can be locally written as
\begin{equation}\label{G3}
  \textbf{x}(v_1,v_2)=(v_1,v_2, \varphi_1(v_1) \varphi_2(v_2))
\end{equation}
or
\begin{equation}\label{G4}
  \textbf{x}(v_1,v_3)=(v_1, \varphi_1(v_1) \varphi_3(v_3),v_3)
\end{equation}
or
\begin{equation}\label{G5}
  \textbf{x}(v_2,v_3)=( \varphi_2(v_2) \varphi_3(v_3),v_2,v_3),
\end{equation}
where $\varphi_i's$ are $C^r,r\geq 1$ smooth functions. The surfaces given by (\ref{G3}), (\ref{G4}) and (\ref{G5}) are called the homothetical surfaces of the first, the second and the third type, respectively.  We have a complete classification result of null Gaussian curvature homothetical surfaces in the following theorem:
\begin{theorem}\label{thm}\cite{Gr2}
Let ${\bf M}$ be a factorable(or homothetical) surface with null Gaussian curvature in $\mathbb{G}_3^1$. If ${\bf M}$ is a factorable surface of the first type (respectively the second type and the third type), then either
\begin{itemize}
\item[(a)] at least one of $\varphi_1$, $\varphi_2$ (respectively $\varphi_1$, $\varphi_3$ and $\varphi_2$, $\varphi_3$) is a constant function, or
\item[(b)] $\varphi_i(v_i)=c_ie^{d_iv_i}$, where $c_i,d_i\in\mathbb{R}\setminus \{0\}$, $i\in \{1,2\}$, (respectively $i\in \{1,3\}$ and $i\in \{2,3\}$) or
\item[(c)] $\varphi_i(v_i)=[(1-m_i)n_i v_i+\lambda_i]^{\frac{1}{1-m_i}}$ where $m_i\neq 0,1,m_i\in \mathbb{R}$ and $m_im_j=1(i\neq j),$ $n_i\in \mathbb{R}\setminus \{0\}$ and $\lambda_i\in \mathbb{R},$ $i\in \{1,2\}$ (respectively $i=\{1,3\}$ and $i\in \{2,3\}$).
\end{itemize}
Conversely, the factorable surfaces satisfying the above cases have null Gaussian curvature.
\end{theorem}
\noindent Since our discussion is about the surfaces with non-degenerate second fundamental form, so in light of this fact, our discussion will be confined to the study of surfaces not falling under the ambit of theorem \ref{thm}.

\noindent Moreover, we see that the first type and the second type homothetical surfaces have up to a sign similar second fundamental form, so we will only discuss first type and  the third type homothetical surfaces \cite{Gr22,Gr2}.
\begin{definition}
  A surface in three-dimensional pseudo-Galilean space $\mathbb{G}_3^1$ is called $\textbf{II}-$harmonic if it satisfies the condition $\Delta^{\textbf{II}}\textbf{x}=\textbf{0}.$\end{definition}
The main results of this paper are:
\begin{theorem}\label{3.3q}
  There are no homothetical surfaces of the first type with non-degenerate second fundamental form satisfying  $\Delta^{\textbf{II}}\textbf{x}_i=\lambda_i \textbf{x}_i,$ where $\lambda_i\in \mathbb{R},$ $i=1,2,3.$ 
\end{theorem}
\begin{theorem}\label{3.3}
  There are no homothetical surfaces of the third type with non-degenerate second fundamental form satisfying  $\Delta^{\textbf{II}}\textbf{x}_i=\lambda_i \textbf{x}_i,$ where $\lambda_i\in \mathbb{R},$ $i=1,2,3.$ 
\end{theorem}
\section{Homothetical surfaces of first type satisfying $\Delta^{II}\textbf{x}_i=\lambda_i\textbf{x}_i$}
 {\bf Proof of theorem \ref{3.3q}:}

Let {\bf x} be a homothetical surface of the first type with non-degenerate second fundamental form in $\mathbb{G}_3^1$ satisfying the condition
\begin{equation}\label{G6}
 \Delta^{II}\textbf{x}_i=\lambda_i\textbf{x}_i,
\end{equation}
where $\lambda_i \in \mathbb{R}, i=1,2,3$ and
\begin{equation*}
  \Delta^{II}\textbf{x}_i=(\Delta^{II}\textbf{x}_1,\Delta^{II}\textbf{x}_2,\Delta^{II}\textbf{x}_3),
\end{equation*}
where $$\textbf{x}_1=v_1, \quad \textbf{x}_2=v_2, \quad \textbf{x}_3= \varphi_1(v_1) \varphi_2(v_2).$$
Now for the homothetical surface given by (\ref{G3}), the coefficients of the second fundamental form are given by
\begin{eqnarray}\label{G7}
 \nonumber&& L=-\frac{\epsilon}{W} \varphi_1^{\prime \prime}\varphi_2, \quad M=-\frac{\epsilon}{W} \varphi_1^{\prime} \varphi_2^{\prime}, \quad N=-\frac{\epsilon}{W}\varphi_1 \varphi_2^{\prime\prime}, \quad \\&& \text{where  } W=\sqrt{\left|1-(\varphi_1  \varphi_2^{\prime})^2\right|}\neq 0.
\end{eqnarray}
The Gaussian curvature $\textbf{K}$ is given by
$$\textbf{K}=\frac{-\epsilon}{W^4}(\varphi_1 \varphi_2 \varphi_1^{\prime\prime} \varphi_2^{\prime\prime}-{ \varphi_1^\prime}^2{ \varphi_2^\prime}^2).$$
Since the surface has non-degenerate second fundamental form everywhere, we have
\begin{equation}\label{Gauss}D=\varphi_1 \varphi_2 \varphi_1^{\prime\prime} \varphi_2^{\prime\prime}-{ \varphi_1^\prime}^2{ \varphi_2^\prime}^2\neq 0,\forall\text{ } v_1,\text{ } v_2\in I.\end{equation}
The Laplacian operator of ${\bf x}_i,(i=1,2,3)$ with the help of (\ref{G2}) turns out to be
\begin{eqnarray}\label{4AG8}
\Delta^{II}{\bf x}&=&
\left(
  \begin{array}{ll}\vspace{.1cm}
    -\frac{W}{\sqrt{D}}\left[\frac{\partial }{\partial v_1}\left(\frac{-\epsilon \varphi_1 \varphi_2^{\prime \prime}}{\sqrt{D}}\right)+\frac{\partial }{\partial v_2}\left(\frac{\epsilon  \varphi_1^\prime  \varphi_2^{\prime}}{\sqrt{D}}\right)\right], \\
   -\frac{W}{\sqrt{D}}\left[\frac{\partial }{\partial v_1}\left(\frac{-\epsilon  \varphi_1^\prime  \varphi_2^{ \prime}}{\sqrt{D}}\right)+\frac{\partial }{\partial v_2}\left(\frac{\epsilon  \varphi_1^{\prime \prime} \varphi_2}{\sqrt{D}}\right)\right], \\
   \varphi_1^\prime \varphi_2 \left(-\frac{W}{\sqrt{D}}\left[\frac{\partial }{\partial v_1}\left(\frac{-\epsilon \varphi_1 \varphi_2^{\prime \prime}}{\sqrt{D}}\right)+\frac{\partial }{\partial v_2}\left(\frac{\epsilon  \varphi_1^\prime  \varphi_2^{\prime}}{\sqrt{D}}\right)\right]\right)+ \\
   \varphi_1 \varphi_2^\prime\left(-\frac{W}{\sqrt{D}}\left[\frac{\partial }{\partial v_1}\left(\frac{-\epsilon  \varphi_1^\prime  \varphi_2^{ \prime}}{\sqrt{D}}\right)+\frac{\partial }{\partial v_2}\left(\epsilon \frac{ \varphi_1^{\prime \prime} \varphi_2}{\sqrt{D}}\right)\right]\right)+2\epsilon \sqrt{1- \varphi_1^2{ \varphi_2^\prime}^2}
  \end{array}
\right).\nonumber\\
\end{eqnarray}
 Since $\textbf{M}$ satisfies (\ref{G6}), equation (\ref{4AG8}) gives rise to the following differential equations
\begin{eqnarray}\label{G9}
    -\frac{W}{\sqrt{D}}\left[\frac{\partial }{\partial v_1}\left(\frac{-\epsilon \varphi_1 \varphi_2^{\prime \prime}}{\sqrt{D}}\right)+\frac{\partial }{\partial v_2}\left( \frac{\epsilon  \varphi_1^\prime  \varphi_2^{\prime}}{\sqrt{D}}\right)\right]=\lambda_1 v_1, \end{eqnarray}
 \begin{eqnarray} \label{G10} -\frac{W}{\sqrt{D}}\left[\frac{\partial }{\partial v_1}\left(\frac{-\epsilon  \varphi_1^\prime  \varphi_2^{ \prime}}{\sqrt{D}}\right)+\frac{\partial }{\partial v_2}\left(\frac{\epsilon  \varphi_1^{\prime \prime} \varphi_2}{\sqrt{D}}\right)\right] =\lambda_2 v_2,\end{eqnarray}
 \begin{eqnarray}\label{G11}\nonumber && \varphi_1^\prime \varphi_2 \left(-\frac{W}{\sqrt{D}}\left[\frac{\partial }{\partial v_1}\left(\frac{-\epsilon \varphi_1 \varphi_2^{\prime \prime}}{\sqrt{D}}\right)+\frac{\partial }{\partial v_2}\left(\frac{\epsilon  \varphi_1^\prime  \varphi_2^{\prime}}{\sqrt{D}}\right)\right]\right)\nonumber \\
  && +\varphi_1 \varphi_2^\prime\left(-\frac{W}{\sqrt{D}}\left[\frac{\partial }{\partial v_1}\left(\frac{-\epsilon  \varphi_1^\prime  \varphi_2^{ \prime}}{\sqrt{D}}\right)+\frac{\partial }{\partial v_2}\left(\frac{\epsilon  \varphi_1^{\prime \prime} \varphi_2}{\sqrt{D}}\right)\right]\right)
\nonumber  \\&&+2\epsilon \sqrt{1- \varphi_1^2{ \varphi_2^\prime}^2}=\lambda_3 \varphi_1 \varphi_2.
\end{eqnarray}
If all the $\lambda_i,(i=1,2,3)$ are distinct, then $\bf M$ is at most of 3-type. By combining (\ref{G9}), (\ref{G10}) and (\ref{G11}), we get
\begin{equation}\label{G12}
   \varphi_1^\prime \varphi_2 \lambda_1 v_1+\varphi_1  \varphi_2^\prime \lambda_2 v_2 +2\epsilon \sqrt{1- \varphi_1^2{ \varphi_2^\prime}^2}= \lambda_3 \varphi_1 \varphi_2.
\end{equation}
Since $ \varphi_1 \varphi_2  \neq 0$, (\ref{G12}) can be written as
\begin{equation}\label{G12a}
  \frac{ \varphi_1^\prime}{\varphi_1}\lambda_1 v_1 + \frac{ \varphi_2^\prime}{\varphi_2}\lambda_2 v_2 + \frac{2\epsilon \sqrt{1-(\varphi_1  \varphi_2^\prime)^2}}{ \varphi_1 \varphi_2 }=\lambda_3.
\end{equation}
According to the choices of constants $\lambda_1,$ $\lambda_2$ and $\lambda_3$, we discuss all the possible cases of $\lambda_i$, $i\in\{1,2,3\}$.\newline
\indent \textbf{Case 1:} Let $\lambda_1=\lambda_2=\lambda_3=0,$ from (\ref{G12a}), we get
$$2\epsilon \sqrt{1- \varphi_1^2{ \varphi_2^\prime}^2}=0, \textrm{ or  } W=0,$$ which is a contradiction to our assumption. Hence there exists no {\bf II}-harmonic homothetical surfaces of first type in $\mathbb{G}_3^1$. \newline
\indent \textbf{Case 2:} Let $\lambda_1=0$, $\lambda_2= 0$ and $\lambda_3\neq 0,$ from (\ref{G12a}), we get
\begin{equation}\label{G13}
 \frac{2\epsilon\sqrt{1- \varphi_1^2{ \varphi_2^\prime}^2}}{ \varphi_1 \varphi_2 }=\lambda_3.
\end{equation}
From (\ref{G13}), we obtain
\begin{equation*}
\frac{4}{ \varphi_1^2}-4 {\varphi_2^\prime}^2-\lambda_3^2  \varphi_2^2=0.
\end{equation*}
Since $\varphi_1$ and $\varphi_2$ are functions of two independent variables, the above equation can be written as
\begin{equation*}
\frac{4}{ \varphi_1^2}=c, \quad 4 {\varphi_2^\prime}^2+\lambda_3^2  \varphi_2^2=c,
\end{equation*}
where $c \in \mathbb{R}\setminus 0$.
Thus, we get
\begin{equation}
 \varphi_1(v_1)=\pm\frac{2}{\sqrt{c}}, \quad  \varphi_2(v_2)=\pm \frac{\sqrt{c}\tan
\left(\frac{1}{2}\lambda_3(v_2+2c_1)\right)}{\lambda_3\sqrt{1+\tan \left(\frac{1}{2}\lambda_3(v_2+2c_1)\right)^2}}.
\end{equation}
In this case the surface may be parameterized as 
\begin{equation}\label{lm0}{\bf x}(v_1,v_2)=\left(v_1,v_2,\left(\pm\frac{2}{\sqrt{c}}\right)\left(\pm \frac{\sqrt{c}\tan\left(\frac{1}{2}\lambda_3(v_2+2c_1)\right)}{\lambda_3\sqrt{1+\tan \left(\frac{1}{2}\lambda_3(v_2+2c_1)\right)^2}}\right)\right).\end{equation}
We observe that the parameterization in (\ref{lm0}) is a contradiction to non-degenerate property as well as to the part $(a)$ of theorem \ref{thm}.

\indent \textbf{Case 3:} Let $\lambda_1=0$, $\lambda_2\neq 0$ and $\lambda_3\neq 0,$ from (\ref{G12a}), we get
\begin{equation}\label{G16}
  \frac{ \varphi_2^\prime}{\varphi_2}\lambda_2 v_2 + \frac{2\epsilon\sqrt{1- \varphi_1^2 { \varphi_2^\prime}^2}}{ \varphi_1 \varphi_2 }=\lambda_3.
\end{equation}
From (\ref{G16}), we obtain
\begin{equation}\label{z1}
\frac{4}{ \varphi_1^2}=c, \quad ( \varphi_2^\prime)^2 +( \varphi_2^\prime \lambda_2 v_2 -\lambda_3  \varphi_2) ^2=c,
\end{equation}
where $c \in \mathbb{R}\setminus 0$. Since the first equation in (\ref{z1}) is constant, so regardless of the second equation of (\ref{z1}), it gives rise to a contradiction to the property of being non-degenerate.  Therefore there exists no parameterization in this case. \newline
\indent \textbf{Case 4:} Let $\lambda_1\neq 0$, $\lambda_2=0$ and $\lambda_3= 0,$ from (\ref{G12a}), we have
\begin{equation*}
  \frac{ \varphi_1^\prime}{\varphi_1}\lambda_1 v_1 + \frac{2\epsilon\sqrt{1- \varphi_1^2 { \varphi_2^\prime}^2}}{ \varphi_1 \varphi_2 }=0.
\end{equation*}
Squaring and adjusting the like terms in above equation, we get
\begin{equation}\label{za1}
( \varphi_1^\prime \lambda_1 v_1)^2=\frac{4}{ \varphi_2^2}-4 \varphi_1^2 \left(\frac{ \varphi_2^\prime}{\varphi_2}\right)^2.
\end{equation}
Differentiating (\ref{za1}) with respect to $v_1$, we get
\begin{equation}
\frac{\lambda_1 ^2}{4}\left(\frac{ \varphi_1^{\prime \prime}}{\varphi_1}v_1^2 +\frac{ \varphi_1^\prime}{\varphi_1}v_1\right)+\left(\frac{ \varphi_2^\prime}{\varphi_2}\right)^2=0.
\end{equation}
Since $\varphi_1$ and $\varphi_2$ are functions of two independent variables, we may write
\begin{equation*}
\frac{\lambda_1 ^2}{4}\left(\frac{ \varphi_1^{\prime \prime}}{\varphi_1}v_1^2 +\frac{ \varphi_1^\prime}{\varphi_1}v_1\right)=-c, \quad \left(\frac{ \varphi_2^\prime}{\varphi_2}\right)^2=c
\end{equation*}
or
\begin{equation}\label{as1}
\left(\frac{ \varphi_1^{\prime \prime}}{\varphi_1}v_1^2 +\frac{ \varphi_1^\prime}{\varphi_1}v_1\right)=\tilde{c}, \quad \left(\frac{ \varphi_2^\prime}{\varphi_2}\right)^2=c,
\end{equation}
where $\tilde{c}=-c\frac{4}{\lambda_1^2}$, $c\in \mathbb{R}$. If $c=0$, then the second equation of (\ref{as1}) implies $\varphi_2=$ constant, which leads to a contradictions. Therefore for $c\in \mathbb{R}\setminus 0,$ we have
\begin{equation}\label{b}
\varphi_1(v_1)=c_1\cos\sqrt{\tilde{c}}\log(v_1)+c_2\sin\sqrt{\tilde{c}}\log(v_1),\quad \varphi_2(v_2)=c_3e^{\pm \sqrt{c}v_2},
\end{equation}
where $c_i \in \mathbb{R}, i\in \{1,2,3\}$.
In this case the, surface may be parameterized as
\begin{equation*}
{\bf x}(v_1,v_2)=\left(v_1,v_2,\left(c_1\cos\sqrt{\tilde{c}}\log(v_1)+c_2\sin\sqrt{\tilde{c}}\log(v_1)\right)\left(c_3e^{\pm \sqrt{c}v_2}\right)\right).
\end{equation*} 
 We observe that the second equation of (\ref{b}) is a contradiction to non-degenerate property with respect to the part $(b)$ of the theorem \ref{thm}. Therefore, there exists no parameterization in this case.
\newline
\indent \textbf{Case 5:} Let $\lambda_1=0$, $\lambda_2\neq0$ and $\lambda_3=0,$ from (\ref{G12a}), we get
\begin{equation*}
  \frac{ \varphi_2^\prime}{\varphi_2}\lambda_2v_2+\frac{2\epsilon\sqrt{1- \varphi_1^2{ \varphi_2^\prime}^2}}{ \varphi_1 \varphi_2 }=0.
\end{equation*}
From above equation, we obtain
\begin{equation}
\frac{4}{ \varphi_1^2}-4( \varphi_2^\prime)^2-( \varphi_2^\prime\lambda_2 v_2)^2=0.
\end{equation}
Since $\varphi_1$ and $\varphi_2$ are functions of independent variables, we can write
\begin{equation}\label{zas1}
(4+\lambda^2_2  v_2^2)( \varphi_2^\prime)^2=c, \quad \frac{4}{ \varphi_1^2}=c,
\end{equation}
where $c \in \mathbb{R}\setminus 0$.
From (\ref{zas1}), we obtain
\begin{equation}
 \varphi_2(v_2)=\pm \frac{\sqrt{c}}{\lambda_2}\sinh ^{-1}\left(\frac{\lambda_2 v_2}{2}\right), \quad  \varphi_1(v_1)=\pm \frac{2}{\sqrt{c}}.
\end{equation}
Therefore the surface may be parameterized as:
\begin{equation}\label{bb}
{\bf x}(v_1,v_2)=\left(v_1,v_2,\left(\pm \frac{2}{\sqrt{c}}\right)\left(\pm \frac{\sqrt{c}}{\lambda_2}\sinh ^{-1}\left(\frac{\lambda_2 v_2}{2}\right)\right)\right).\end{equation}
It is clearly visible that the parameterization in (\ref{bb}) gives rise to a similar type of contradiction as in case 2.

\indent  \textbf{Case 6:} Let $\lambda_1\neq0$, $\lambda_2\neq0$ and $\lambda_3=0,$ from (\ref{G12a}), we get
\begin{equation}\label{qa1}
  \frac{ \varphi_1^\prime}{\varphi_1}\lambda_1 v_1+\frac{ \varphi_2^\prime}{\varphi_2}\lambda_2 v_2+\frac{2\epsilon \sqrt{1- \varphi_1^2{ \varphi_2^\prime}^2}}{ \varphi_1 \varphi_2 }=0.
\end{equation}
There is no suitable solution of (\ref{qa1}). Hence there exists no parameterization in this case also.

\indent  \textbf{Case 7:} Let $\lambda_1\neq0$, $\lambda_2=0$ and $\lambda_3 \neq 0,$ from (\ref{G12a}), we have
\begin{equation}\label{,1}
  \frac{ \varphi_1^\prime}{\varphi_1}\lambda_1 v_1+\frac{2\epsilon \sqrt{1- \varphi_1^2{ \varphi_2^\prime}^2}}{ \varphi_1 \varphi_2 }=\lambda_3.
\end{equation}
Differentiating (\ref{,1}) with respect to $v_2$, we get
\begin{equation}\label{,2}
\frac{2\epsilon \sqrt{1- \varphi_1^2{ \varphi_2^\prime}^2}}{ \varphi_1 \varphi_2 }=c,
\end{equation}
where $c\in \mathbb{R}$. If $c=0$, then it is a contradiction to $W\neq 0$. So assuming $c\neq 0$, squaring and adjusting the like terms in (\ref{,2}), we get
\begin{equation}\label{,3}
\frac{ 1}{\left(\varphi_1\right)^2}=\varphi_2^2 \left(c+\left(\frac{\varphi_2^\prime}{\varphi_2^2}\right)^2\right).
\end{equation}
Since $\varphi_1$ and $\varphi_2$ are functions of two independent variables, we can write
\begin{equation}\label{,4}
\frac{ 1}{\left(\varphi_1\right)^2}=c_1, \textrm{ }\varphi_2^2 \left(c+\left(\frac{\varphi_2^\prime}{\varphi_2^2}\right)^2\right)=c_1,
\end{equation}
where $c_1 \in \mathbb{R}\setminus 0$. In any way, regardless of the solution of second equation of (\ref{,4}), the first equation in  (\ref{,4}) will contradict to the non-degenerate property. Hence there exists no parameterization in this case.

\indent  \textbf{Case 8:} Let $\lambda_1\neq0$, $\lambda_2\neq 0$ and $\lambda_3 \neq 0,$ from (\ref{G12a}), we have
\begin{equation}\label{,5}
  \frac{ \varphi_1^\prime}{\varphi_1}\lambda_1 v_1+
  \frac{ \varphi_2^\prime}{\varphi_2}\lambda_2 v_2+\frac{2\epsilon \sqrt{1- \varphi_1^2{ \varphi_2^\prime}^2}}{ \varphi_1 \varphi_2 }=\lambda_3.
\end{equation}
Differentiating (\ref{,5}) with respect to $v_1$ and $v_2$, we get
\begin{equation}
\frac{2\epsilon \sqrt{1- \varphi_1^2{ \varphi_2^\prime}^2}}{ \varphi_1 \varphi_2 }=c
\end{equation}
If $c=0$, then it is a contradiction to $W\neq0$. Suppose $c\neq 0$, from (\ref{,5}), we obtain
 \begin{equation}\label{,6}
  \frac{ \varphi_1^\prime}{\varphi_1}\lambda_1 v_1+
  \frac{ \varphi_2^\prime}{\varphi_2}\lambda_2 v_2+c=\lambda_3.
\end{equation}
Since $\varphi_1$ and $\varphi_2$ are functions of two independent variables, we can write 
\begin{equation}\label{,7}
 \frac{ \varphi_1^\prime}{\varphi_1}\lambda_1 v_1=c_1,\textrm{ }
 \lambda_3-\frac{ \varphi_2^\prime}{\varphi_2}\lambda_2 v_2+c=c_1.
\end{equation}
We can easily see that the solutions of (\ref{,7}) does not satisfy (\ref{,5}). 
 This completes the proof of the theorem \ref{3.3q}.
\section{Homothetical surfaces of third type satisfying $\Delta^{II}\textbf{x}_i=\lambda_i\textbf{x}_i$}
{\bf Proof of theorem \ref{3.3}:}

Let {\bf x} be a homothetical surface of the third type with non-degenerate second fundamental form in $\mathbb{G}_3^1$ satisfying the condition
\begin{equation}\label{G22}
 \Delta^{II}\textbf{x}_i=\lambda_i\textbf{x}_i,
\end{equation}
where $\lambda_i \in \mathbb{R}, i=1,2,3$ and
\begin{equation*}
  \Delta^{II}\textbf{x}_i=(\Delta^{II}\textbf{x}_1,\Delta^{II}\textbf{x}_2,\Delta^{II}\textbf{x}_3),
\end{equation*}
where $$\textbf{x}_1= \varphi_1(v_2) \varphi_2(v_3), \quad \textbf{x}_2=v_2, \quad \textbf{x}_3=v_3.$$
For the homothetical surface given by (\ref{G5}), the coefficients of the second fundamental form are given by
\begin{eqnarray}\label{G23}
\nonumber  &&L=-\frac{\epsilon}{\mathcal{W}} \varphi_1^{\prime \prime}\varphi_2, \quad M=-\frac{\epsilon}{\mathcal{W}} \varphi_1^{\prime} \varphi_2^{\prime}, \quad N=-\frac{\epsilon}{\mathcal{W}} \varphi_1 \varphi_2^{\prime\prime}, \\
  && \text{where  } \mathcal{W}=\sqrt{\left|( \varphi_1^\prime  \varphi_2) ^2-(\varphi_1  \varphi_2^{\prime})^2\right|}\neq 0.
\end{eqnarray}
The Gaussian curvature $\textbf{K}$ is given by
$$\textbf{K}=\frac{-\epsilon}{\mathcal{W}^4}(\varphi_1 \varphi_2 \varphi_1^{\prime\prime} \varphi_2^{\prime\prime}-{ \varphi_1^\prime}^2{ \varphi_2^\prime}^2).$$
Since the surface has non-degenerate second fundamental form, we have
$$ \mathcal{D}=\varphi_1 \varphi_2 \varphi_1^{\prime\prime} \varphi_2^{\prime\prime}-{ \varphi_1^\prime}^2{ \varphi_2^\prime}^2\neq 0,\forall\text{ } v_2, v_3\in I.$$
In this case, the Laplacian operator of ${\bf x}_i,i=1,2,3$ with the help of (\ref{G2}) turns out to be
\begin{eqnarray}\label{4AG24}
\Delta^{II}{\bf x}&=&
\left(
  \begin{array}{ll}\vspace{.1cm}
   \varphi_1^\prime \varphi_2\left(-\frac{\mathcal{W}}{\sqrt{ \mathcal{D}}}\left[\frac{\partial }{\partial v_2}\left(\frac{-\epsilon  \varphi_1 \varphi_2^{\prime \prime}}{\sqrt{ \mathcal{D}}}\right)+\frac{\partial }{\partial v_3}\left(\frac{\epsilon  \varphi_1^\prime  \varphi_2^{\prime}}{\sqrt{ \mathcal{D}}}\right)\right]\right) \\
    +\varphi_1 \varphi_2^\prime\left(-\frac{\mathcal{W}}{\sqrt{ \mathcal{D}}}\left[\frac{\partial }{\partial v_2}\left(\frac{-\epsilon  \varphi_1^\prime  \varphi_2^{ \prime}}{\sqrt{ \mathcal{D}}}\right)+\frac{\partial }{\partial v_3}\left(\frac{\epsilon  \varphi_1^{\prime \prime} \varphi_2}{\sqrt{ \mathcal{D}}}\right)\right]\right)\\+2\epsilon \sqrt{{ \varphi_1^\prime}^2 \varphi_2^2- \varphi_1^2{ \varphi_2^\prime}^2},\\
    -\frac{\mathcal{W}}{\sqrt{ \mathcal{D}}}\left[\frac{\partial }{\partial v_2}\left(\frac{-\epsilon  \varphi_1 \varphi_2^{\prime \prime}}{\sqrt{ \mathcal{D}}}\right)+\frac{\partial }{\partial v_3}\left(\frac{\epsilon  \varphi_1^\prime  \varphi_2^{\prime}}{\sqrt{ \mathcal{D}}}\right)\right],\vspace{.1cm} \\ 
   -\frac{\mathcal{W}}{\sqrt{ \mathcal{D}}}\left[\frac{\partial }{\partial v_2}\left(\frac{-\epsilon  \varphi_1^\prime  \varphi_2^{ \prime}}{\sqrt{ \mathcal{D}}}\right)+\frac{\partial }{\partial v_3}\left(\frac{\epsilon  \varphi_1^{\prime \prime} \varphi_2}{\sqrt{ \mathcal{D}}}\right)\right],
  \end{array}
\right).\nonumber\\
&&
\end{eqnarray}
Since $\textbf{M}$ satisfies (\ref{G22}), equation (\ref{4AG24}) gives rise to the following differential equations
\begin{eqnarray}\label{G27}\nonumber && \varphi_1^\prime \varphi_2\left(-\frac{\mathcal{W}}{\sqrt{ \mathcal{D}}}\left[\frac{\partial }{\partial v_2}\left(\frac{-\epsilon  \varphi_1 \varphi_2^{\prime \prime}}{\sqrt{ \mathcal{D}}}\right)+\frac{\partial }{\partial v_3}\left(\frac{\epsilon  \varphi_1^\prime  \varphi_2^{\prime}}{\sqrt{ \mathcal{D}}}\right)\right]\right) \nonumber\\
  && + \varphi_1 \varphi_2^\prime\left(-\frac{\mathcal{W}}{\sqrt{ \mathcal{D}}}\left[\frac{\partial }{\partial v_2}\left(\frac{-\epsilon  \varphi_1^\prime  \varphi_2^{ \prime}}{\sqrt{ \mathcal{D}}}\right)+\frac{\partial }{\partial v_3}\left(\frac{\epsilon  \varphi_1^{\prime \prime} \varphi_2}{\sqrt{ \mathcal{D}}}\right)\right]\right)\\
  &&+2\epsilon \sqrt{{ \varphi_1^\prime}^2 \varphi_2^2- \varphi_1^2{ \varphi_2^\prime}^2}=\lambda_1 \varphi_1 \varphi_2 ,\nonumber 
\end{eqnarray}
\begin{eqnarray}\label{G25}
    -\frac{\mathcal{W}}{\sqrt{ \mathcal{D}}}\left[\frac{\partial }{\partial v_2}\left(\frac{-\epsilon  \varphi_1 \varphi_2^{\prime \prime}}{\sqrt{ \mathcal{D}}}\right)+\frac{\partial }{\partial v_3}\left(\frac{\epsilon  \varphi_1^\prime  \varphi_2^{\prime}}{\sqrt{ \mathcal{D}}}\right)\right]=\lambda_2 v_2, \end{eqnarray}
    \begin{eqnarray} \label{G26} -\frac{\mathcal{W}}{\sqrt{ \mathcal{D}}}\left[\frac{\partial }{\partial v_2}\left(\frac{-\epsilon  \varphi_1^\prime  \varphi_2^{ \prime}}{\sqrt{ \mathcal{D}}}\right)+\frac{\partial }{\partial v_3}\left(\frac{\epsilon  \varphi_1^{\prime \prime} \varphi_2}{\sqrt{ \mathcal{D}}}\right)\right] =\lambda_3 v_3.\end{eqnarray}
Depending upon the nature of $\lambda_i,(i=1,2,3)$, $\textbf{M}$ is at most of 3-type. By combining (\ref{G27}), (\ref{G25}) and (\ref{G26}), we get
\begin{equation*}
   \varphi_1^\prime \varphi_2 \lambda_2 v_2+ \varphi_1  \varphi_2^\prime \lambda_3 v_3 +2\epsilon \sqrt{{ \varphi_1^\prime}^2 \varphi_2^2- \varphi_1^2{ \varphi_2^\prime}^2}= \lambda_1  \varphi_1 \varphi_2 .
\end{equation*}
The above equation can be written as
\begin{equation}\label{G28}
  \frac{ \varphi_1^\prime}{\varphi_1} \lambda_2 v_2+ \frac{ \varphi_2^\prime}{\varphi_2} \lambda_3 v_3 +\frac{2\epsilon\sqrt{{ \varphi_1^\prime}^2 \varphi_2^2- \varphi_1^2{ \varphi_2^\prime}^2}}{ \varphi_1 \varphi_2 }= \lambda_1.
\end{equation}
\indent  \textbf{Case 1:} Let $\lambda_2=0$, $\lambda_3=0$ and $\lambda_1=0,$ from (\ref{G28}), we get
\begin{equation*}
  \frac{2\epsilon\sqrt{{ \varphi_1^\prime}^2 \varphi_2^2- \varphi_1^2{ \varphi_2^\prime}^2}}{ \varphi_1 \varphi_2 }=0,
\end{equation*}
which is a contradiction to the non-vanishing assumption of $\mathcal{W}$. Hence there exists no ${\bf II}$-harmonic homothetical surfaces of the third type in $\mathbb{G}_3^1.$ \newline
\indent  \textbf{Case 2:} Let $\lambda_2=0$, $\lambda_3=0$ and $\lambda_1\neq 0,$ from (\ref{G28}), we get
\begin{equation*}
  \frac{2\epsilon\sqrt{{ \varphi_1^\prime}^2 \varphi_2^2- \varphi_1^2{ \varphi_2^\prime}^2}}{ \varphi_1 \varphi_2 }=\lambda_1.
\end{equation*}
From above equation, we derive
\begin{equation}\label{ha}
\left\{\left(\frac{ \varphi_1^\prime}{\varphi_1}\right)^2 -\frac{\lambda_1}{8}\right\}-\left\{\left(\frac{ \varphi_2^\prime}{\varphi_2}\right)^2 +\frac{\lambda_1}{8}\right\}=0.
\end{equation}
Since $\varphi_1$ and $\varphi_2$ are functions of two independent variables, equaiton (\ref{ha}) can be written as
\begin{equation}
\left(\frac{ \varphi_1^\prime}{\varphi_1}\right)^2 -\frac{\lambda_1}{8}=c, \quad 
\left(\frac{ \varphi_2^\prime}{\varphi_2}\right)^2 +\frac{\lambda_1}{8}=c,
\end{equation}
where $c\in \mathbb{R}$.
Hence we get
\begin{eqnarray}
 \label{fa}\varphi_1(v_2)&=&c_1e^{\pm\frac{1}{2\sqrt{2}}\sqrt{\lambda_1 +8c} v_2}\\
\label{fb} \varphi_2(v_3)&=&c_1e^{\pm\frac{1}{2\sqrt{2}}\sqrt{\lambda_1-8c} v_3},\quad c<\frac{\lambda_1}{8},
\end{eqnarray}
where $c_1\in \mathbb{R}$. 

\indent \textbf{Case 3:} Let $\lambda_2=0$, $\lambda_3\neq 0$ and $\lambda_1\neq 0,$ from (\ref{G28}), we get
\begin{equation*}
 \frac{ \varphi_2^\prime}{\varphi_2} \lambda_3 v_3 +\frac{2\epsilon\sqrt{{ \varphi_1^\prime}^2 \varphi_2^2- \varphi_1^2{ \varphi_2^\prime}^2}}{ \varphi_1 \varphi_2 }= \lambda_1.
\end{equation*}
From above equation, we obtain
\begin{equation*}
4\left\{\left(\frac{ \varphi_1^\prime}{\varphi_1}\right)^2-\left(\frac{ \varphi_2^\prime}{\varphi_2}\right)^2\right\}=\left\{\lambda_1 -\frac{ \varphi_2^\prime}{\varphi_2}\lambda_3 v_3\right\}^2
\end{equation*}
or
\begin{equation}\label{bg1}
4\left(\frac{ \varphi_1^\prime}{\varphi_1}\right)^2-\lambda_1^2-4\left(\frac{ \varphi_2^\prime}{\varphi_2}\right)^2-\left(\frac{ \varphi_2^\prime}{\varphi_2}\lambda_3 v_3\right)^2+2\lambda_3\lambda_1 v_3 \left(\frac{ \varphi_2^\prime}{\varphi_2}\right)=0.
\end{equation}
Since $\varphi_1$ and $\varphi_2$ are functions of two independent variables, we can write (\ref{bg1}) as
\begin{equation*}
4\left(\frac{ \varphi_1^\prime}{\varphi_1}\right)^2-\lambda_1^2=c, \quad
 4\left(\frac{ \varphi_2^\prime}{\varphi_2}\right)^2+\left(\frac{ \varphi_2^\prime}{\varphi_2}\lambda_3 v_3\right)^2-2\lambda_3\lambda_1 v_3 \left(\frac{ \varphi_2^\prime}{\varphi_2}\right)=c,
\end{equation*}where $c\in \mathbb{R}$.
Thus we have
\begin{eqnarray}
 \label{fc}\varphi_1(v_2)&=&c_1 e^{\pm \frac{1}{2}\sqrt{\lambda_1^2 +c}v_2},\\ 
 \label{fd} \varphi_2(v_3)&=&\frac{1}{\lambda_3^2}(4+\lambda_3^2v_3^2)^{\frac{\lambda_1}{2\lambda_3}}\left(n^2v_3+nm\right)^{-\frac{n}{\lambda_3^2}}e^{\pm\lambda_3\lambda_1 \tanh^{-1}\left(\frac{\lambda_3\lambda_1 v_3}{m}\right)+c_2}.\nonumber \\
\end{eqnarray}
where
\begin{equation*}
m=\sqrt{4c+(\lambda_3\lambda_1^2+\lambda_3^2\lambda_1)v_3^2}, \quad n=\sqrt{\lambda_3\lambda_1^2+\lambda_3^2\lambda_1}\textrm{ and }c_1,c_2 \in \mathbb{R}.
\end{equation*} 

\indent \textbf{Case 4:} Let $\lambda_2\neq0$, $\lambda_3= 0$ and $\lambda_1= 0,$ from (\ref{G28}), we get
\begin{equation*}
  \frac{ \varphi_1^\prime}{\varphi_1} \lambda_2 v_2 +\frac{2\epsilon\sqrt{{ \varphi_1^\prime}^2 \varphi_2^2- \varphi_1^2{ \varphi_2^\prime}^2}}{ \varphi_1 \varphi_2 }= 0.
\end{equation*}
The above equation can be rewritten as 
\begin{equation}\label{df8}
\left(\frac{ \varphi_1^\prime}{\varphi_1}\right)^2 \left(4-\lambda_2^2 v_2^2\right)-4\left(\frac{ \varphi_2^\prime}{\varphi_2}\right)^2=0.
\end{equation}
We can write (\ref{df8}) in the following form:
\begin{equation*}
\left(\frac{ \varphi_1^\prime}{\varphi_1}\right)^2 \left(4-\lambda_2^2 v_2^2\right)=c, \quad 4\left(\frac{ \varphi_2^\prime}{\varphi_2}\right)^2=c,
\end{equation*}
where $c\in \mathbb{R}$. If $c=0$, then from the second part of above equation we obtain $\varphi_2=$ constant, which leads to a contradiction. Thus, for $c\in \mathbb{R}\setminus 0$, we have
\begin{equation}\label{fe}
 \varphi_1(v_2)=\frac{c_1}{\lambda_2}e^{\pm \sqrt{c}\sin^{-1}\left(\frac{\lambda_2 v_2}{2}\right)}, \quad  \varphi_2(v_3)=c_2e^{\pm \frac{\sqrt{c}}{2}v_3}.
\end{equation}
where $c_1,c_2$ are non-zero constants.

\indent \textbf{Case 5:} Let $\lambda_2=0$, $\lambda_3\neq 0$ and $\lambda_1= 0,$ from (\ref{G28}), we get
\begin{equation*}
  \frac{ \varphi_2^\prime}{\varphi_2} \lambda_3v_3  +\frac{2\epsilon\sqrt{{ \varphi_1^\prime}^2 \varphi_2^2- \varphi_1^2{ \varphi_2^\prime}^2}}{ \varphi_1 \varphi_2 }= 0.
\end{equation*}
On the similar lines as in case 4, we can easily obtain
\begin{equation}\label{fg}
 \varphi_1(v_2)=c_1e^{\pm \frac{\sqrt{c}}{2}v_2}, \quad  \varphi_2(v_3)=\frac{c_2}{\lambda_3}e^{\pm \sqrt{c}\sinh^{-1}\left(\frac{\lambda_3 v_3}{2}\right)}.
\end{equation}
where $c_1,c_2 \in \mathbb{R}\setminus 0$. From (\ref{fa}) and (\ref{fb}), (\ref{fc}) and (\ref{fd}), (\ref{fe}), (\ref{fg}), we see that there exists at least one $\varphi_i$, $i \in \{1,2\}$ of the similar form as in part $(b)$ of the theorem \ref{thm}, leading to a contradiction to the non-degenerate property. Therefore there exists no required parameterization from case 2 to case 5.

\indent\textbf{Case 6:} Let $\lambda_2\neq0$, $\lambda_3\neq 0$ and $\lambda_1= 0,$ from (\ref{G28}), we get
\begin{equation}\label{fda}
\frac{ \varphi_1^\prime}{\varphi_1} \lambda_2 v_2+\frac{ \varphi_2^\prime}{\varphi_2} \lambda_3 v_3 +\frac{2\epsilon\sqrt{{ \varphi_1^\prime}^2 \varphi_2^2- \varphi_1^2{ \varphi_2^\prime}^2}}{ \varphi_1 \varphi_2 }= 0.\end{equation}
There exists no non-trivial analytic solution of (\ref{fda}) other than $\varphi_1=\varphi_2=$ constant, but this assumption again causes a contradiction. 

\indent\textbf{Case 7:} Let $\lambda_2\neq0$, $\lambda_3= 0$ and $\lambda_1\neq 0,$ from (\ref{G28}), we get
\begin{equation}\label{,a}
\frac{ \varphi_1^\prime}{\varphi_1} \lambda_2 v_2 +\frac{2\epsilon\sqrt{{ \varphi_1^\prime}^2 \varphi_2^2- \varphi_1^2{ \varphi_2^\prime}^2}}{ \varphi_1 \varphi_2 }=\lambda_3.\end{equation}
Squaring and adjusting the like terms, we have
\begin{equation*}
\left(\frac{ \varphi_1^\prime}{\varphi_1}\right)^2-\left(\frac{ \varphi_2^\prime}{\varphi_2}\right)^2= \frac{1}{4}\left(\lambda_3-\frac{ \varphi_1^\prime}{\varphi_1} \lambda_2 v_2 \right)^2\end{equation*}
or
\begin{equation*}
\left(\frac{ \varphi_1^\prime}{\varphi_1}\right)^2-
\frac{1}{4}\left(\lambda_3-\frac{ \varphi_1^\prime}{\varphi_1} \lambda_2 v_2 \right)^2=\left(\frac{ \varphi_2^\prime}{\varphi_2}\right)^2.\end{equation*}
The above equation can be written as
\begin{equation}\label{,b}
\left(\frac{ \varphi_2^\prime}{\varphi_2}\right)^2=c,\textrm{ }
\left(\frac{ \varphi_1^\prime}{\varphi_1}\right)^2-
\frac{1}{4}\left(\lambda_3-\frac{ \varphi_1^\prime}{\varphi_1} \lambda_2 v_2 \right)^2=c,
\end{equation}
where $c\in \mathbb{R}$. If $c=0$, from the first equation of (\ref{,b}), $\varphi_2=$ constant leads to a contradiction. Suppose $c\neq0$, from the first equation of (\ref{,b}), we obtain $\varphi_2=e^{\pm \sqrt{c}v_3}$ which is again a contradiction to part $(b)$ of theorem \ref{thm}.

\indent\textbf{Case 8:} Let $\lambda_2\neq0$, $\lambda_3\neq 0$ and $\lambda_1 \neq0,$ from (\ref{G28}), we get
\begin{equation}\label{fda}
\frac{ \varphi_1^\prime}{\varphi_1} \lambda_2 v_2+\frac{ \varphi_2^\prime}{\varphi_2} \lambda_3 v_3 +\frac{2\epsilon\sqrt{{ \varphi_1^\prime}^2 \varphi_2^2- \varphi_1^2{ \varphi_2^\prime}^2}}{ \varphi_1 \varphi_2 }= \lambda_1.\end{equation}
Following the similar steps as in case $8$ of (\ref{,5}), we arrive at similar types of contradictions. Therefore, there exists no parameterization in this case also.
This completes the proof of the theorem \ref{3.3}.

{\bf Acknowledgement:} I am very thankful to the anonymous referees for their valuable comments and suggestions which helped a lot to improve the quality of this paper.

\end{document}